\documentclass{amsart}
\usepackage{amsmath,amssymb}
\usepackage{yhmath}

\newcommand{\bdis}{\begin{displaymath}}
\newcommand{\edis}{\end{displaymath}}
\newcommand{\be}{\begin{equation}}
\newcommand{\ee}{\end{equation}}
\newcommand{\mbb}{\mathbb}
\newcommand{\mcal}{\mathcal}

\newcommand{\vp}{\varphi}

\newcommand{\zf}{\zeta\left(\frac{1}{2}+it\right)}
 
\newcommand{\zfn}{\zeta\left(\frac{1}{2}+it_\nu\right)}   
\newcommand{\zfnn}{\zeta\left(\frac{1}{2}+it_{\nu+1}\right)}

\theoremstyle{definition}

\theoremstyle{remark}
\newtheorem{remark}[]{Remark}

\newtheorem*{mydef11}{{\bf Theorem 1}}

\newtheorem*{mydef12}{{\bf Theorem 2}}

\newtheorem*{mydef15}{{\bf Theorem 5}} 

\newtheorem*{mydef16}{{\bf Theorem 6}}

\newtheorem*{mydef41}{{\bf Corollary 1}}

\newtheorem*{mydef51}{{\bf Lemma 1}}

\newtheorem*{mydef52}{{\bf Lemma 2}}

\newtheorem*{mydef53}{{\bf Lemma 3}}

\newtheorem*{mydef54}{{\bf Lemma 4}}

\newtheorem*{mydef57}{{\bf Lemma 7}}

\newtheorem*{mydefI}{\bf Corporate consequence}

\numberwithin{equation}{section}

\begin{document}

\title[Jacob's ladders, Hardy-Littlewood integral (1918) \dots ]{Jacob's ladders, Hardy-Littlewood integral (1918) and new asymptotic functional equations for Euler's Gamma function together with the tenth equivalent of the Fermat-Wiles theorem}

\author{Jan Moser}

\address{Department of Mathematical Analysis and Numerical Mathematics, Comenius University, Mlynska Dolina M105, 842 48 Bratislava, SLOVAKIA}

\email{jan.mozer@fmph.uniba.sk}


\begin{abstract}
In this paper new $\Gamma$-functional is constructed upon the basis of the set of almost linear increments of the Hardy-Littlewood integral. This functional generates a $\Gamma$-equivalent of the Fermat-Wiles theorem and also new set of factorization formulae for Euler's $\Gamma$-function.   
\end{abstract}
\maketitle

\section{Introduction} 

\subsection{} 

Let us remind that we have obtained new class of the increments of the Hardy-Littlewood integral\footnote{See \cite{1}.}
\be \label{1.1} 
J(T)=\int_0^T \left|\zf\right|^2{\rm d}t 
\ee 
in our paper \cite{6}. Namely, it is true, that for every sufficiently big $T>0$ and for every fixed $k\in\mbb{N}$, there is a sequence 
\be \label{1.2} 
\{\overset{r}{T}(T)\}_{r=1}^k:\ T<\overset{1}{T}(T)<\overset{2}{T}(T)<\dots<\overset{k}{T}(T)
\ee 
such, that the following formula 
\be \label{1.3} 
\int_{\overset{r-1}{T}}^{\overset{r}{T}}\left|\zf\right|^2{\rm d}t\sim (1-c)\overset{r-1}{T},\ T\to\infty 
\ee 
holds true, where $c$ is the Euler's constant. 

\begin{remark}
This almost linear formula (w.r.t the variable $\overset{r-1}{T}$) represents the new property of the Hardy-Littlewood integral (\ref{1.1}) on the continuum set of segments 
\be \label{1.4}  
\begin{split}
& [\overset{r-1}{T},\overset{r}{T}]:\ \overset{r}{T}-\overset{r-1}{T}\sim (1-c)\pi(\overset{r}{T}), \\ 
& T\to\infty,\ r=1,\dots,k
\end{split}
\ee  
(see \cite{6}, (3.3)), where 
\be \label{1.5} 
\pi(\overset{r}{T})\sim\frac{\overset{r}{T}}{\ln\overset{r}{T}},\ T\to\infty, 
\ee  
and the symbol $\pi(x)$ stands for the prime-counting function. 
\end{remark}

Next, in the papers \cite{7} -- \cite{9} we have defined, on the basis of the formula (\ref{1.3}), some functionals that generate nine equivalents of the Fermat-Wiles theorem. 

\subsection{} 

Let us remind the Euler's $\Gamma$-function 
\be \label{1.6} 
\Gamma(x)=\int_0^\infty t^{x-1}e^{-t}{\rm d}t,\ x>0. 
\ee  
In this paper we have constructed the functional 
\be \label{1.7} 
\lim_{\tau\to \infty}\ln\left\{
\frac{\Gamma([\frac{x}{1-c}\tau]^1)}{\Gamma(\frac{x}{1-c}\tau)}
\right\}^{\frac{1}{\tau}}=x, 
\ee  
for every $x>0$, where 
\be \label{1.8} 
\left[\frac{x}{1-c}\tau\right]^1=\vp_{1}^{-1}\left( \frac{x}{1-c}\tau\right). 
\ee  
This formula associates any fixed ray 
\be \label{1.9} 
y(\tau;x)=\frac{x}{1-c}\tau 
\ee 
with the slope 
\be \label{1.10} 
0<\arctan\frac{x}{1-c}<\frac{\pi}{2}
\ee 
with the unique element $x>0$.  

\subsection{} 

First, the following condition follows from (\ref{1.7}) 
\be \label{1.11} 
\lim_{\tau\to\infty}
\left\{
\frac{\Gamma([\frac{x^n+y^n}{z^n}\frac{\tau}{1-c}]^1)}{\Gamma(\frac{x^n+y^n}{z^n}\frac{\tau}{1-c})}
\right\}^{\frac{1}{\tau}}\not= e
\ee 
on the class of all Fermat's rationals 
\be \label{1.12} 
\frac{x^n+y^n}{z^n},\ x,y,z\in\mbb{N},\ n\geq 3, 
\ee 
and this condition represents the $\Gamma$-equivalent of the Fermat-Wiles theorem. 

\begin{remark}
The $\Gamma$-equivalent of the Fermat-Wiles theorem (\ref{1.11}) contains Euler's notion of $\Gamma(t)$ and $c$ as well as the set of all Fermat rationals. These notions are connected by means of the reverse iterations $\vp_1^{-1}(t)$ of the Jacob's ladder. 
\end{remark} 

\subsection{} 

Next, we have proved the following asymptotic functional equations of new kind 
\be \label{1.13} 
\Gamma(\tau)\left\{
\prod_{\tau\leq n\leq\overset{1}{\tau}}e^{d(n)}
\right\}\sim \Gamma(\overset{1}{\tau}), 
\ee 
\be \label{1.14} 
\begin{split}
& \{\Gamma(\tau)\}^{\frac{1+c}{\pi}}
\left\{
\prod_{\tau<t_\nu\leq \overset{1}{\tau}}\exp\left[\zfn\zfnn\right]
\right\}\sim \{\Gamma(\overset{1}{\tau})\}^{\frac{1+c}{\pi}}, 
\end{split}
\ee 
\be \label{1.15} 
\{\Gamma(\tau)\}^{\frac{1}{\pi}}
\left\{
\prod_{\tau<t_\nu\leq \overset{1}{\tau}}\exp\left[\zfn\right]
\right\}\sim \{\Gamma(\overset{1}{\tau})\}^{\frac{1}{\pi}}, 
\ee  
where 
\bdis 
\overset{1}{\tau}=[\tau]^1=\vp_1^{-1}(\tau),\ \tau\to\infty, 
\edis  
by means of the functional (\ref{1.7}). The symbol $d(n)$ stands for the number of divisors of $n$ and 
\be \label{1.16} 
\{t_\nu\}_{\nu=1}^\infty 
\ee  
denotes the Gram's sequence. 

\begin{remark}
The formula (\ref{1.13}) can be viewed as, for example, 
\begin{itemize}
	\item[(a)] the asymptotic factorization formula for the function 
	\be \label{1.17} 
	\Gamma(\vp_1^{-1}(\tau)), 
	\ee 
	\item[(b)] the asymptotic transformation 
	\be \label{1.18} 
	\Gamma(\tau)\xrightarrow{\Phi_1}\Gamma(\overset{1}{\tau}). 
	\ee 
\end{itemize} 
The same is true also in the cases of formulae (\ref{1.14}) and (\ref{1.15}), namely 
\be\label{1.19} 
\{\Gamma(\tau)\}^{\frac{1+c}{\pi}}\xrightarrow{\Phi_2}\{\Gamma(\overset{1}{\tau})\}^{\frac{1+c}{\pi}}, 
\ee 
\be\label{1.20} 
\{\Gamma(\tau)\}^{\frac{1}{\pi}}\xrightarrow{\Phi_3}\{\Gamma(\overset{1}{\tau})\}^{\frac{1}{\pi}}. 
\ee
\end{remark} 

\begin{remark}
One can characterize the transformation (\ref{1.18})\footnote{As well as (\ref{1.19}) and (\ref{1.20}).} by the following way: External multiplication of the heterogeneous elements on the left-hand side of eq. (\ref{1.13}) generates the internal transformation 
\bdis 
\tau\to \overset{1}{\tau} 
\edis 
of the argument of Euler's Gamma-function. 
\end{remark} 

\subsection{} 

In this paper we use the following notions of our works \cite{2} -- \cite{5}: 
\begin{itemize}
	\item[{\tt (a)}] Jacob's ladder $\vp_1(T)$, 
	\item[{\tt (b)}] direct iterations of Jacob's ladders 
	\bdis 
	\begin{split}
		& \vp_1^0(t)=t,\ \vp_1^1(t)=\vp_1(t),\ \vp_1^2(t)=\vp_1(\vp_1(t)),\dots , \\ 
		& \vp_1^k(t)=\vp_1(\vp_1^{k-1}(t))
	\end{split}
	\edis 
	for every fixed natural number $k$, 
	\item[{\tt (c)}] reverse iterations of Jacob's ladders 
	\be \label{1.21} 
	\begin{split}
		& \vp_1^{-1}(T)=\overset{1}{T},\ \vp_1^{-2}(T)=\vp_1^{-1}(\overset{1}{T})=\overset{2}{T},\dots, \\ 
		& \vp_1^{-r}(T)=\vp_1^{-1}(\overset{r-1}{T})=\overset{r}{T},\ r=1,\dots,k, 
	\end{split} 
	\ee 
	where, for example, 
	\be \label{1.22} 
	\vp_1(\overset{r}{T})=\overset{r-1}{T}
	\ee  
	for every fixed $k\in\mbb{N}$ and every sufficiently big $T>0$. Next, we use also the properties of reverse iterations: 
	\be \label{1.23}
	\begin{split} 
		& \overset{r}{T}-\overset{r-1}{T}\sim(1-c)\pi(\overset{r}{T});\ \pi(\overset{r}{T})\sim\frac{\overset{r}{T}}{\ln \overset{r}{T}},\ r=1,\dots,k,\ T\to\infty, \\ 
		& \overset{0}{T}=T<\overset{1}{T}(T)<\overset{2}{T}(T)<\dots<\overset{k}{T}(T), \\ 
		& T\sim \overset{1}{T}\sim \overset{2}{T}\sim \dots\sim \overset{k}{T},\ T\to\infty. 
	\end{split}
	\ee 
\end{itemize} 

\begin{remark}
	The asymptotic behaviour of the points 
	\bdis 
	\{T,\overset{1}{T},\dots,\overset{k}{T}\}
	\edis  
	is as follows: at $T\to\infty$ these points recede unboundedly each from other and all together are receding to infinity. Hence, the set of these points behaves at $T\to\infty$ as one-dimensional Friedmann-Hubble expanding Universe. 
\end{remark} 

\section{List of equivalents of the Fermat-Wiles theorem} 

\subsection{} 

Our paper \cite{7} contains the following $\zeta$-equivalents of the Fermat-Wiles theorem: 
\be \label{2.1} 
\lim_{\tau\to\infty}\frac{1}{\tau}\int_{\frac{x^n+y^n}{z^n}\frac{\tau}{1-c}}^{\left[\frac{x^n+y^n}{z^n}\frac{\tau}{1-c}\right]^1}\left|\zf\right|^2{\rm d}t\not= 1, 
\ee   
\be \label{2.2} 
\lim_{\tau\to\infty}
\frac
{\int_{(x^n+y^n)\frac{\tau}{1-c}}^{\left[(x^n+y^n)\frac{\tau}{1-c}\right]^1}\left|\zf\right|^2{\rm d}t}
{\int_{z^n\frac{\tau}{1-c}}^{\left[z^n\frac{\tau}{1-c}\right]^1}\left|\zf\right|^2{\rm d}t}\not= 1, 
\ee 
\be \label{2.3} 
\lim_{\tau\to\infty}\ln\left\{
\int_{\exp(\frac{x^n+y^n}{z^n}\ln \tau)}^{[\exp(\frac{x^n+y^n}{z^n}\ln \tau)]^1}\left|\zf\right|^2{\rm d}t
\right\}^{\frac{1}{\ln\tau}}\not=1,  
\ee  
\be \label{2.4} 
\lim_{\tau\to\infty}
\frac
{\ln\{\int_{\exp((x^n+y^n)\ln\tau)}^{[\exp((x^n+y^n)\ln\tau)]^1}|\zf|^2{\rm d}t\}}
{\ln\{\int_{\exp(z^n\ln\tau)}^{[\exp(z^n\ln\tau)]^1}|\zf|^2{\rm d}t\}} 
\not=1. 
\ee  

\subsection{} 

Let us remind the Dirichlet $D(x)$-function 
\be \label{2.5} 
D(x)=\sum_{n\leq x}d(n), 
\ee 
where $d(n)$ is the number of divisors of the positive integer $n$ and  
\be \label{2.6} 
D(x)=D(N),\ x\in [N,N+1),\ \forall N\in\mbb{N}. 
\ee 
Our paper \cite{8} contains two $D$-equivalents of the Fermat-Wiles theorem:  
\be \label{2.7} 
\lim_{\tau\to\infty}\frac{1}{\tau}\left\{
D\left(\left[\frac{x^n+y^n}{z^n}\frac{\tau}{1-c}\right]^1\right)-D\left(\frac{x^n+y^n}{z^n}\frac{\tau}{1-c}\right)
\right\}\not=1, 
\ee 
\be \label{2.8} 
\lim_{\tau\to\infty}\ln\left\{
D\left(\left[\exp\left(\frac{x^n+y^n}{z^n}\ln\tau\right)\right]^1\right)
-D\left(\exp\left(\frac{x^n+y^n}{z^n}\ln\tau\right)\right)
\right\}^{\frac{1}{\ln\tau}}\not=1, 
\ee  
and also one $\zeta$-equivalent of the Fermat-Wilson theorem, namely: 
\be \label{2.9} 
\begin{split}
	& \lim_{\rho\to\infty}\frac{1}{\rho}\int_{\vp_1(\frac{x^n+y^n}{z^n}\frac{\rho}{(1-c)\sigma(l)})}^{\frac{x^n+y^n}{z^n}\frac{\rho}{(1-c)\sigma(l)}}\left\{
	\sigma(l)\left|\zf\right|^2+(1-c)|S_1(t)|^{2l}
	\right\}{\rm d}t\not=1,  
\end{split}
\ee  
where 
\be \label{2.10}  
S_1(t)=\frac{1}{\pi}\int_0^t\arg\zf{\rm d}t 
\ee 
and $\sigma(l)$ denotes the Selberg's constants. 

\subsection{} 

Next, let us remind the Titchmarsh's functions 
\be \label{2.11} 
\mcal{T}_1(X)=\sum_{t_\nu\leq X} \zfn, 
\ee  
and 
\be \label{2.12} 
\mcal{T}_2(X)=\sum_{t_\nu\leq X} \zfn\zfnn, 
\ee 
where 
\be \label{2.13} 
\{t_\nu\}_{\nu=1}^\infty,\ t_{\nu+1}-t_\nu\sim\frac{2\pi}{\ln t_\nu},\ \nu\to\infty 
\ee 
is the Gram's sequence. Our paper \cite{9} contains further two ($\mcal{T}_1$ and $\mcal{T}_2$) equivalents of the Fermat-Wiles theorem: 
\be \label{2.14} 
\lim_{\tau\to\infty}\frac{1}{\tau}\left\{\mcal{T}_1\left(\left[\frac{x^n+y^n}{z^n}\frac{\tau}{1-c}\right]^1\right)-\mcal{T}_1\left(\frac{x^n+y^n}{z^n}\frac{\tau}{1-c}\right)\right\}\not=\frac{1}{\pi}, 
\ee  
\be \label{2.15} 
\lim_{\tau\to\infty}\frac{1}{\tau}\left\{\mcal{T}_2\left(\left[\frac{x^n+y^n}{z^n}\frac{\tau}{1-c}\right]^1\right)-\mcal{T}_2\left(\frac{x^n+y^n}{z^n}\frac{\tau}{1-c}\right)\right\}\not=\frac{1+c}{\pi}. 
\ee 

\subsection{} 

We give the tenth equivalent of the Fermat-Wiles theorem in the present paper, namely the following Gamma-equivalent 
\be \label{2.16} 
\lim_{\tau\to\infty}\ln\left\{
\frac{\Gamma([\frac{x^n+y^n}{z^n}\frac{\tau}{1-c}]^1)}{\Gamma(\frac{x^n+y^n}{z^n}\frac{\tau}{1-c})} 
\right\}^{\frac{1}{\tau}}\not=1.
\ee  
It is assumed, of course, that every of the conditions (\ref{2.1}) -- (\ref{2.4}), (\ref{2.7}) -- (\ref{2.9}) and (\ref{2.14}) -- (\ref{2.16}) holds true on the class of all Fermat's rationals.\footnote{Compare the subsection 1.3.}

\section{The first lemma as the asymptotic Newton-Leibniz formula for the function $|\zf|^2$}

\subsection{} 

Let us remind the Stirling formula 
\be \label{3.1} 
\Gamma(x)=x^{x-\frac 12}e^{-x}\sqrt{2\pi}e^{\frac{\Theta}{12x}},\ \Theta\in (0,1) 
\ee 
for the Euler's Gamma-function. In the logarithmic form we have 
\bdis 
\ln\Gamma(x)=\left(x-\frac 12\right)\ln x-x+\frac 12\ln 2\pi+\mcal{O}\left(\frac 1x\right). 
\edis 
For our purpose we can use simplified version of previous formula, namely 
\be \label{3.2} 
\ln\Gamma(x)=x\ln x-x+\mcal{O}(\ln x),\ x\to\infty. 
\ee  

\subsection{} 

Let us remind our formula\footnote{See \cite{2} and \cite{3}.} 
\be \label{3.3} 
\begin{split}
& \int_0^T\left|\zf\right|^2{\rm d}t=\vp_1(T)\ln\{\vp_1(T)\}+ \\ 
& (c-\ln 2\pi)\vp_1(T)+c_0+\mcal{O}\left(\frac{\ln T}{T}\right) 
\end{split}
\ee 
giving us the infinite set of the almost exact representations of the Hardy-Littlewood integral (\ref{1.1}). 

Now, the substitution 
\be \label{3.4} 
T\to\overset{r}{T},\ r=1,\dots,k
\ee 
for eery fixed $k\in\mbb{N}$ in (\ref{3.3}) results in\footnote{See (\ref{1.22}).} 
\be \label{3.5} 
\int_0^{\overset{r}{T}}\left|\zf\right|^2{\rm d}t=\overset{r-1}{T}\ln\overset{r-1}{T}+(c-\ln 2\pi)\overset{r-1}{T}+\mcal{O}(1), 
\ee 
while the substitution 
\be \label{3.6} 
x\to\overset{r-1}{T}
\ee 
in (\ref{3.2}) gives\footnote{See (\ref{1.23}).}
\be \label{3.7} 
\ln\Gamma(\overset{r-1}{T})=\overset{r-1}{T}\ln\overset{r-1}{T}-\overset{r-1}{T}+\mcal{O}(\ln T). 
\ee 
Subtracting (\ref{3.7}) from (\ref{3.5}) we have 
\be \label{3.8} 
\int_0^{\overset{r}{T}}\left|\zf\right|^2{\rm d}t-\ln\Gamma(\overset{r-1}{T})=(1+c-\ln 2\pi)\overset{r-1}{T}+\mcal{O}(\ln T), 
\ee  
and translation $r\to r+1$ in (\ref{3.8}) results in 
\be \label{3.9} 
\int_0^{\overset{r+1}{T}}\left|\zf\right|^2{\rm d}t-\ln\Gamma(\overset{r}{T})=(1+c-\ln 2\pi)\overset{r}{T}+\mcal{O}(\ln T). 
\ee  
And consequently, by sbtraction of (\ref{3.8}) from (\ref{3.9}), we obtain 
\be \label{3.10} 
\begin{split}
& \int_{\overset{r}{T}}^{\overset{r+1}{T}}\left|\zf\right|^2{\rm d}t-\ln\Gamma(\overset{r}{T})+\ln\Gamma(\overset{r-1}{T})= \\ 
& (1+c-\ln 2\pi)(\overset{r}{T}-\overset{r-1}{T})+\mcal{O}(\ln T). 
\end{split} 
\ee 

\subsection{} 

Let us remind the following: 
\begin{itemize}
	\item[(a)] our almost linear formula\footnote{See \cite{6}, (3.4), (3.6); $0<\delta$ is sufficiently small.}
	\be \label{3.11} 
	\int_{\overset{r-1}{T}}^{\overset{r}{T}}\left|\zf\right|^2{\rm d}t=(1-c)\overset{r-1}{T}+\mcal{O}(T^{\frac 13 + \delta}), 
	\ee  
	that implies 
	\be \label{3.12} 
	\begin{split}
	& \int_{\overset{r}{T}}^{\overset{r+1}{T}}\left|\zf\right|^2{\rm d}t=\int_{\overset{r-1}{T}}^{\overset{r}{T}}\left|\zf\right|^2{\rm d}t+ \\ 
	& (1-c)(\overset{r}{T}-\overset{r-1}{T})+\mcal{O}(T^{\frac 13 + \delta}); 
	\end{split}
	\ee 
	\item[(b)] and also our next formula\footnote{See \cite{2}, comp. \cite{3}, (7.2).} 
	\be \label{3.13} 
	T-\vp_1(T)\sim (1-c)\pi(T),\ T\to\infty, 
	\ee  
	that implies, see (\ref{1.23}) 
	\be \label{3.14} 
	\overset{r}{T}-\overset{r-1}{T}\sim (1-c)\frac{\overset{r}{T}}{\ln \overset{r}{T}}\sim (1-c)\frac{T}{\ln T}=\mcal{O}\left(\frac{T}{\ln T}\right). 
	\ee
\end{itemize} 

Finally, we get from the eq. (\ref{3.10}) (see also (\ref{3.12}) and (\ref{3.14})) the following lemma. 

\begin{mydef51}
\be \label{3.15} 
\begin{split}
& \ln\frac{\Gamma(\overset{r}{T})}{\Gamma(\overset{r-1}{T})}=\int_{\overset{r-1}{T}}^{\overset{r}{T}}\left|\zf\right|^2{\rm d}t+\mcal{O}\left(\frac{T}{\ln T}\right), \\ 
& r=1,\dots,k,\ T\to\infty. 
\end{split}
\ee 
\end{mydef51} 

\begin{remark}
Rewriting the formula (\ref{3.15}) into the form 
\be\label{3.16} 
\ln\Gamma(\overset{r}{T})-\ln\Gamma(\overset{r-1}{T})\sim \int_{\overset{r-1}{T}}^{\overset{r}{T}}\left|\zf\right|^2{\rm d}t
\ee 
we exhibit, that the eq. (\ref{3.15}) is an asymptotic Newton-Leibniz formula for the function 
\bdis 
\left|\zf\right|^2
\edis  
on the continuum set of segments\footnote{See (\ref{1.23}).} 
\bdis 
[\overset{r-1}{T}(T),\overset{r}{T}(T)]:\ \overset{r}{T}(T)-\overset{r-1}{T}(T)\sim (1-c)\pi[\overset{r}{T}(T)],\ T\to\infty. 
\edis 
\end{remark} 

\section{Next lemmas and therorem 1} 

\subsection{} 

In the following text we shall use our formula (\ref{3.15}) with $r=1$, for example, i. e. 
\be \label{4.1} 
\begin{split}
	& \ln\frac{\Gamma(\overset{1}{T})}{\Gamma(T)}=\int_{T}^{\overset{1}{T}}\left|\zf\right|^2{\rm d}t+\mcal{O}\left(\frac{T}{\ln T}\right), \\ 
	& T=\overset{0}{T},\ T>T_0>0, 
\end{split}
\ee 
where $T_0$ is sufficiently big and 
\be \label{4.2} 
\overset{1}{T}=[T]^1=\vp_1^{-1}(T). 
\ee 
By making use of the substitution 
\be \label{4.3} 
T=\frac{x}{1-c}\tau,\ \tau\in \left(\frac{1-c}{x}T_0,+\infty\right),\ x>0 
\ee 
in eq. (\ref{4.1}) we obtain the following statement. 

\begin{mydef52}
\be \label{4.4} 
\begin{split}
& \ln\frac{\Gamma([\frac{x}{1-c}\tau]^1)}{\Gamma(\frac{x}{1-c}\tau)}= 
\int_{\frac{x}{1-c}\tau}^{[\frac{x}{1-c}\tau]^1}\left|\zf\right|^2{\rm d}t+\mcal{O}\left(\frac{\tau}{\ln\tau}\right), \\ 
& \tau\in)(\tau_1(x),+\infty),\ \tau_1(x)=\max\left\{\left(\frac{1-c}{x}\right)^2,(T_0)^2\right\}, 
\end{split}
\ee 
for every fixed $x>0$. 
\end{mydef52} 

Since\footnote{See \cite{7}, (4.6).} 
\be \label{4.5} 
\lim_{\tau\to\infty}\frac{1}{\tau}\int_{\frac{x}{1-c}\tau}^{[\frac{x}{1-c}\tau]^1}\left|\zf\right|^2{\rm d}t=x, 
\ee 
then it follows from (\ref{4.4}), that the next lemma holds true.  

\begin{mydef53}
\be \label{4.6} 
\lim_{\tau\to\infty}\ln\left\{\ln\frac{\Gamma([\frac{x}{1-c}\tau]^1)}{\Gamma(\frac{x}{1-c}\tau)}\right\}^{\frac{1}{\tau}}=x 
\ee 
for every fixed $x>0$, where 
\be \label{4.7} 
\left[\frac{x}{1-c}\tau\right]^1=\vp_1^{-1}\left(\frac{x}{1-c}\tau\right). 
\ee 
\end{mydef53} 

\subsection{} 

Now, by making use of the substitution 
\be \label{4.8} 
x\to \frac{x^n+y^n}{z^n} 
\ee 
in eq. (\ref{4.6}), we obtain the following. 

\begin{mydef54}
The formula 
\be \label{4.9} 
\lim_{\tau\to\infty}\ln\left\{
\frac{\Gamma([\frac{x^n+y^n}{z^n}\frac{\tau}{1-c}]^1)}{\Gamma(\frac{x^n+y^n}{z^n}\frac{\tau}{1-c})} 
\right\}^{\frac{1}{\tau}}=\frac{x^n+y^n}{z^n} 
\ee  
for every fixed Fermat's rational 
\bdis 
\frac{x^n+y^n}{z^n} . 
\edis 
\end{mydef54} 

Finally, the Theorem 1 holds true. 

\begin{mydef11}
The $\Gamma$-condition 
\be \label{4.10} 
\lim_{\tau\to\infty}\ln\left\{
\frac{\Gamma([\frac{x^n+y^n}{z^n}\frac{\tau}{1-c}]^1)}{\Gamma(\frac{x^n+y^n}{z^n}\frac{\tau}{1-c})} 
\right\}^{\frac{1}{\tau}}\not=1
\ee 
on the class of all Fermat's rationals represents the $\Gamma$-equivalent of the Fermat-Wiles theorem. 
\end{mydef11} 

It follows from (\ref{4.10}) immediately, that 

\begin{mydef41}
The condition 
\be \label{4.11} 
\lim_{\tau\to\infty}\left\{
\frac{\Gamma([\frac{x^n+y^n}{z^n}\frac{\tau}{1-c}]^1)}{\Gamma(\frac{x^n+y^n}{z^n}\frac{\tau}{1-c})} 
\right\}^{\frac{1}{\tau}}\not=e
\ee 
on the class of all Fermat's rationals represents the second form of the $\Gamma$-equivalent of the Fermat-Wiles theorem. 
\end{mydef41}

\begin{remark}
The second form of the $\Gamma$-equivalent (\ref{4.11}) contains the Euler's notions $\Gamma(t)$ as well as the number $e$ together with the set of Fermat's rationals and all of these are connected by means of the reverse iteration $\vp_1^{-1}(t)$ of the Jacob's ladder. 
\end{remark} 

\section{New factorization formulas for the Euler's Gamma-function} 

\subsection{} 

Let us start with the formula\footnote{See \cite{8}, (5.10).} 
\be \label{5.1} 
\lim_{\tau\to\infty}\frac{1}{\tau}\left\{
D\left(\left[\frac{x}{1-c}\tau\right]^1\right)-D\left(\frac{x}{1-c}\tau\right)
\right\}=x,\ x>0. 
\ee  
Of course, it is sufficient to study the simplest case $x=1-c$, that is 
\be \label{5.2} 
\lim_{\tau\to\infty}\frac{1}{\tau}\{D(\overset{1}{\tau})-D(\tau)\}=1-c. 
\ee  
Next, we have ((\ref{4.6}), $x=1-c$) 
\be \label{5.3} 
\lim_{\tau\to\infty}\ln\left\{
\frac{\Gamma(\overset{1}{\tau})}{\Gamma(\tau)}
\right\}^{\frac{1}{\tau}}=1-c, 
\ee 
and, of course, 
\be \label{5.4} 
\lim_{\tau\to\infty}\frac{D(\overset{1}{\tau})-D(\tau)}{\ln \frac{\Gamma(\overset{1}{\tau})}{\Gamma(\tau)}}=1. 
\ee 
Thus 
\be \label{5.5} 
D(\overset{1}{\tau})-D(\tau)\sim \ln\frac{\Gamma(\overset{1}{\tau})}{\Gamma(\tau)},\ \tau\to\infty, 
\ee  
and\footnote{See (\ref{2.5}).} 
\be \label{5.6} 
\sum_{\tau\leq n\leq\overset{1}{\tau}}d(n)\sim \ln\frac{\Gamma(\overset{1}{\tau})}{\Gamma(\tau)}. 
\ee 
Consequently, the following theorem holds true. 

\begin{mydef12}
\be \label{5.7} 
\frac{\Gamma(\overset{1}{\tau})}{\Gamma(\tau)}\sim \prod_{\tau\leq n\leq\overset{1}{\tau}}e^{d(n)},\ \tau\to\infty. 
\ee 
\end{mydef12} 

\subsection{} 

Next, we use the formula\footnote{See \cite{9}, (5.14).} 
\be \label{5.8} 
\lim_{\tau\to\infty}\frac{1}{\tau}
\left\{
\mcal{T}_2\left(\left[\frac{x}{1-c}\tau\right]^1\right)-\mcal{T}_2\left(\frac{x}{1-c}\tau\right)
\right\}=\frac{1+c}{\pi}x
\ee 
as basic one. In the case $x=1-c$ we obtain 
\be \label{5.9} 
\lim_{\tau\to\infty}\frac{1}{\tau}\{\mcal{T}_2(\overset{1}{\tau})-\mcal{T}_2(\tau)\}=\frac{1+c}{\pi}(1-c). 
\ee  
This formula together with the eq. (\ref{5.3}) imply 
\be \label{5.10} 
\mcal{T}_2(\overset{1}{\tau})-\mcal{T}_2(\tau)\sim \ln\left\{\frac{\Gamma(\overset{1}{\tau})}{\Gamma(\tau)}\right\}^{\frac{1+c}{\pi}},\ \tau\to\infty, 
\ee 
where\footnote{See (\ref{2.12}).}  
\be \label{5.11} 
\mcal{T}_2(\overset{1}{\tau})-\mcal{T}_2(\tau)=\sum_{\tau<t_\nu\leq \overset{1}{\tau}}\zfn\zfnn,\ \tau\to\infty. 
\ee 
As a consequence we have the following theorem. 

\begin{mydef53}
\be \label{5.12} 
\left\{\frac{\Gamma(\overset{1}{\tau})}{\Gamma(\tau)}\right\}^{\frac{1+c}{\pi}}\sim \prod_{\tau\leq t_\nu\leq\overset{1}{\tau}}\exp\left[\zfn\zfnn\right],\ \tau\to\infty. 
\ee 
\end{mydef53} 

\subsection{} 

And finally, we use the formula\footnote{See \cite{9}, (4.12).} 
\be \label{5.13} 
\lim_{\tau\to\infty}\frac{1}{\tau}
\left\{
\mcal{T}_1\left(\left[\frac{x}{1-c}\tau\right]^1\right)-\mcal{T}_1\left(\frac{x}{1-c}\tau\right)
\right\}=\frac{1}{\pi}x
\ee 
and its consequence for $x=1-c$ 
\be \label{5.14} 
\lim_{\tau\to\infty}\frac{1}{\tau}
\left\{
\mcal{T}_1(\overset{1}{\tau})-\mcal{T}_1(\tau)
\right\}=\frac{1-c}{\pi}, 
\ee  
that together with (\ref{5.3}) imply the following statement. 

\begin{mydef54}
\be \label{5.15} 
\left\{\frac{\Gamma(\overset{1}{\tau})}{\Gamma(\tau)}\right\}^{\frac{1}{\pi}}\sim \prod_{\tau\leq t_\nu\leq \overset{1}{\tau}}\exp\left[\zfn\right],\ \tau\to\infty. 
\ee 
\end{mydef54} 

\begin{remark}
The ordering of the formulae (\ref{5.1}), (\ref{5.8}) and (\ref{5.13}) is defined by the inequalities 
\be \label{5.16} 
1>\frac{1+c}{\pi}>\frac{1}{\pi}. 
\ee 
\end{remark} 

\begin{remark}
It is useful to remind the following two asymptotic formulae 
\be \label{5.17} 
\overset{1}{\tau}-\tau\sim (1-c)\pi(\overset{1}{\tau}),\ \tau\to\infty, 
\ee 
and 
\be \label{5.18} 
t_{\nu+1}-t_\nu\sim \frac{2\pi}{\ln t_\nu},\ \nu\to\infty. 
\ee 
in the connection with (\ref{5.7}), (\ref{5.12}) and (\ref{5.15}). 
\end{remark} 

\section{Further formulas} 

\subsection{} 

Since\footnote{See (\ref{5.12}) and (\ref{5.15}).} 
\be \label{6.1} 
\frac{\Gamma(\overset{1}{\tau})}{\Gamma(\tau)}\sim\prod_{\tau< t_\nu\leq \overset{1}{\tau}}\exp\left[\frac{\pi}{1+c}\zfn\zfnn\right], 
\ee  
\be \label{6.2} 
\frac{\Gamma(\overset{1}{\tau})}{\Gamma(\tau)}\sim\prod_{\tau< t_\nu\leq \overset{1}{\tau}}\exp\left[\pi\zfn\right], 
\ee  
then it follows from our formulae (\ref{5.7}). (\ref{6.1}) and (\ref{6.2}). 
\begin{mydefI}
\be \label{6.3}  
\begin{split}
& \prod_{\tau<n \leq \overset{1}{\tau}}e^{d(n)}\sim \prod_{\tau< t_\nu\leq \overset{1}{\tau}}\exp\left[\frac{\pi}{1+c}\zfn\zfnn\right]\sim \\ 
& \prod_{\tau< t_\nu\leq \overset{1}{\tau}}\exp\left[\pi\zfn\right],\ \tau\to\infty. 
\end{split}
\ee 
\end{mydefI} 

\subsection{} 

As long as our formula (\ref{6.2}) holds true for every sufficiently big $\tau$, then this one holds true also for 
\bdis 
\tau\to\overset{1}{\tau}, \overset{2}{\tau}, \dots,\overset{k}{\tau};\ [\overset{1}{\tau}]^1=\overset{2}{\tau},\dots 
\edis 
in other words, we have 
\be \label{6.4} 
\frac{\Gamma(\overset{r}{\tau})}{\Gamma(\overset{r-1}{\tau})}\sim \prod_{\overset{r-1}{\tau}< t_\nu\leq \overset{r}{\tau}}\exp\left[\pi\zfn\right]. 
\ee 
Now we can multiply all the equations in (\ref{6.4}) to obtain 

\begin{mydef41} 
\be \label{6.5} 
\Gamma(\overset{k}{\tau})\sim\Gamma(\tau)\left\{\prod_{\tau< t_\nu\leq \overset{k}{\tau}}\exp\left[\pi\zfn\right]
\right\}
\ee 
holds true and therefore we have 
\be \label{6.6} 
\begin{split}
& \Gamma(\overset{s}{\tau})\sim\Gamma(\overset{r}{\tau})\left\{\prod_{\overset{r}{\tau}< t_\nu\leq \overset{s}{\tau}}\exp\left[\pi\zfn\right]\right\}, \\ 
& 1\leq r\leq s\leq k 
\end{split}
\ee 
for every fixed $k\in\mbb{N}$ and $\tau\to\infty$.  
\end{mydef41} 

\begin{remark}
Of course, formulae corresponding to (\ref{6.5}) and (\ref{6.6}) hold true also in the cases (\ref{5.7}) and (\ref{6.1}). 
\end{remark} 

\subsection{} 

Next, following the Euler's formula 
\be \label{6.7} 
\Gamma(x+1)=x\Gamma(x),\ x>0, 
\ee  
we obtain 
\be \label{6.8} 
\frac{\Gamma(\overset{r}{\tau}+1)}{\Gamma(\overset{r}{\tau})}-\frac{\Gamma(\overset{r-1}{\tau}+1)}{\Gamma(\overset{r-1}{\tau})}=\overset{r}{\tau}-\overset{r-1}{\tau}, \ r=1,\dots,k 
\ee 
and, consequently,  the formula 
\be \label{6.9} 
\overset{r}{\tau}-\overset{r-1}{\tau}\sim (1-c)\pi(\overset{r}{\tau}),\ \tau\to\infty
\ee 
gives the result 
\be \label{6.10} 
\pi(\overset{r}{\tau})\sim\frac{1}{1-c}\left\{
\frac{\Gamma(\overset{r}{\tau}+1)}{\Gamma(\overset{r}{\tau})}-\frac{\Gamma(\overset{r-1}{\tau}+1)}{\Gamma(\overset{r-1}{\tau})}
\right\},\ r=1,\dots,k,\ \tau\to\infty. 
\ee 
Now, summing up all the formulas (\ref{6.10}) together with\footnote{See (\ref{1.23}).}
\be \label{6.11} 
\pi(\overset{r}{\tau})\sim\pi(\tau),\ \tau\to\infty
\ee 
gives the statement 
\begin{mydef15}
\be \label{6.12} 
\pi(\tau)\sim \frac{1}{(1-c)k}\left\{
\frac{\Gamma(\overset{k}{\tau}+1)}{\Gamma(\overset{k}{\tau})}-\frac{\Gamma(\tau+1)}{\Gamma(\tau)},\ \tau\to\infty 
\right\}
\ee 
for every fixed $k\in\mbb{N}$. 
\end{mydef15} 

\subsection{} 

Making another choice from the set of formulas \{(\ref{5.7}), (\ref{6.1}), (\ref{6.2})\} gives us new results. Namely, if we take 
\be \label{6.13} 
\begin{split}
& \Gamma([\tau+1]^1)\sim \Gamma(\tau)\left\{
\prod_{\tau+1<t_\nu\leq[\tau+1]^1}\exp\left[\pi\zfn\right]
\right\}, \\ 
& \Gamma([\tau]^1)\sim \Gamma(\tau)\left\{
\prod_{\tau<t_\nu\leq[\tau]^1}\exp\left[\pi\zfn\right]
\right\}, 
\end{split}
\ee  
together with the Euler's formula (\ref{6.7}), then we obtain 

\begin{mydef16} 
\be \label{6.14} 
\frac{\Gamma([\tau+1]^1)}{\Gamma([\tau]^1)}\sim\tau\frac{\prod_{\tau+1<t_\nu\leq[\tau+1]^1}\exp\left[\pi\zfn\right]}{\prod_{\tau<t_\nu\leq[\tau]^1}\exp\left[\pi\zfn\right]},\ \tau\to\infty. 
\ee 
\end{mydef16} 

\begin{remark}
In the case of decomposition 
\be \label{6.15} 
\{\tau<t_\nu\leq[\tau]^1\}=\{\tau<t_\nu\leq \tau+1\}\cup \{\tau+1<t_\nu\leq \tau\}
\ee 
it is true, that the number of segments 
\bdis 
[t_\nu,t_{\nu+1}]\subset(\tau,\tau+1] 
\edis 
is equal\footnote{See (\ref{2.13}).} to 
\be\label{6.16} 
\sim\frac{1}{2\pi}\ln\tau \to\infty \ \mbox{as}\ \tau\to\infty. 
\ee 
\end{remark} 

\begin{remark}
Further, the Legendre's formula 
\be \label{6.17} 
\frac{\Gamma(2\tau)}{\Gamma(\tau)\Gamma(\tau+\frac 12)}=\frac{1}{\sqrt{\pi}}2^{2\tau-1}
\ee  
generates, together with (\ref{6.2}), the following factorization 
\be \label{6.18} 
\begin{split}
& \frac{\Gamma([2\tau]^1)}{\Gamma([\tau]^1)\Gamma([\tau+\frac 12]^1)}\sim\frac{1}{\sqrt{\pi}}2^{2\tau+1}\left\{
\prod_{2\tau<t_\nu\leq [2\tau]^1}\exp\left[\pi\zfn\right]
\right\}\times \\ 
& \left\{
\prod_{\tau<t_\nu\leq [\tau]^1}\exp\left[\pi\zfn\right]
\right\}^{-1}\times \\ 
& \left\{
\prod_{\tau+\frac 12<t_\nu\leq [\tau+\frac 12]^1}\exp\left[\pi\zfn\right]
\right\}^{-1}, 
\end{split}
\ee  
and so on. 
\end{remark} 

\section{On a localized $\zeta$-equivalent of the Fermat-Wiles theorem} 

\subsection{} 

Let us remind our previous results: 
\begin{itemize}
	\item[(A)] By \cite{7}, Lemma 4 we have 
	\be \label{7.1} 
	\lim_{\tau\to\infty}\frac{1}{\tau}\int_{\frac{x^n+y^n}{z^n}\frac{\tau}{1-c}}^{\left[\frac{x^n+y^n}{z^n}\frac{\tau}{1-c}\right]^1}\left|\zf\right|^2{\rm d}t=\frac{x^n+y^n}{z^n} 
	\ee  
	for every fixed Fermat's rational. 
	\item[(B)] By \cite{7}, Theorem 4 we have the $\zeta$-condition 
	\be \label{7.2} 
	\lim_{\tau\to\infty}\frac{1}{\tau}\int_{\frac{x^n+y^n}{z^n}\frac{\tau}{1-c}}^{\left[\frac{x^n+y^n}{z^n}\frac{\tau}{1-c}\right]^1}\left|\zf\right|^2{\rm d}t\not= 1
	\ee 
	on the class of all Fermat's rationals that is the $\zeta$-equivalent of the Fermat-Wiles theorem. 
\end{itemize} 

\subsection{} 

Since (\ref{7.2}), we have for every small positive and fixed $\epsilon$ and every 
\be \label{7.4} 
\frac{x^n+y^n}{z^n}:\ \left|\frac{x^n+y^n}{z^n}-1\right|\geq\epsilon
\ee 
the result: 

\begin{mydef57}
The $\zeta$-condition 
\be \label{7.5} 
\lim_{\tau\to\infty}\frac{1}{\tau}\int_{\frac{x^n+y^n}{z^n}\frac{\tau}{1-c}}^{\left[\frac{x^n+y^n}{z^n}\frac{\tau}{1-c}\right]^1}\left|\zf\right|^2{\rm d}t\not= 1
\ee 
for 
\be\label{7.6} 
\frac{x^n+y^n}{z^n}\in (1-\epsilon,1+\epsilon) 
\ee  
represents the localized $\zeta$-equivalent of the Fermat-Wiles theorem. 
\end{mydef57} 

\begin{remark}
Of course, there are similar localizations for all other equivalents of the Fermat-Wiles theorem, which we have obtained in our papers \cite{7} -- \cite{9} and also in the present paper. 
\end{remark}

I would like to thank Michal Demetrian for his moral support of my study of Jacob's ladders.

\end{document}